\documentclass[11pt]{article}
\usepackage{amsfonts}
\usepackage{amsmath}

\newtheorem{lemma}{Lemma}
\begin{document}
\begin{titlepage}
\title{Cosh-Gordon equation and quasi-Fuchsian groups}  
\author{V.V. Fock}
\date{}

\end{titlepage}
\maketitle
\abstract{We show propositions in favour of relations between the phase space of 3D gravity, moduli of quasi-Fuchsian groups, global solutions of cosh-Gordon equations and minimal surfaces in hyperbolic spaces.}
\section{Introduction}
The main motivation of the paper was an attempt to study the phase space of Einstein 3D gravity with negative cosmological constant.

The paper does not aim anything else but just a few observations, computations and propositions relating the problem to quasi-Fuchsian groups, affine Lie groups, hyper-K\"aler structures and integrable systems. 

A large number of published papers is devoted to every subject listed above, so we do not feel competent enough to single out a few. We are going to describe here a number of possible approaches to the study of the phase space one of which is very close to the that of C.Taubes \cite{Taubes}.

Let us first define the main actors:

\paragraph{Kleinian and quasi-Fuchsian groups.} (See e.g. \cite{Maskit} for more details.) Recall that a finitely generated discrete torsion free subgroup of the group $PSL(2,\mathbb R)$ is called Fuchsian. Since $PSL(2,\mathbb R)$ is the symmetry group of the hyperbolic plane $\mathbb H^2$ the space of Fuchsian groups naturally decomposes into connected components corresponding to the topology of the quotient of $\mathbb H^2$. The connected component of the space of Fuchsian groups considered up to conjugation and corresponding to a surface $\Sigma$ (not nescesarily closed) is called the Teichm\"uller space of $\Sigma$ and is denoted by $\mathcal T(\Sigma)$. Of course this space can also be interpretated as the space of complex structures on $\Sigma$ or as the space of hyperbolic metrics on $\Sigma$, both considered up to the action of the group $\mathit{Diff}^0$ of diffeomorpisms of $\Sigma$ homotopy equivalent to the identity. The Teichm\"uller space $\mathcal T(\Sigma)$ is a subset of the moduli space $\mathcal M(\Sigma,PSL(2,\mathbb R))$ of flat $PSL(2,\mathbb R)$ connections on the surface $\Sigma$.

Analogously, a Kleinian group is a finitely generated discrete torsion free subgroup of  the group $PSL(2, \mathbb C)$. The space of Kleinian groups considered up to a conjugation also has connected components distinguished by the topology of the quotients of the three-dimensional hyperbolic space $\mathbb H^3$. It is denoted by $\mathcal K(X)$ where $X$ is a three-dimensional manifold. This space can be also considered as the space of complete hyperbolic metrics on $X$ up to diffeomorphisms homotopy equivalent to the identity. Obviously any finitely generated subgroup of a Kleinian group is Kleinian.

If a quotient of the hyperbolic plane $\mathbb H^2$ is a Riemann surface $\Sigma$ then the quotient of the hyperbolic space $\mathbb H^3$ by the same group is topologically a product $\Sigma \times \mathbb R$. Thus a Fuchsian group can be considered as a Kleinian one and therefore the Teichm\"uller space $\mathcal T(\Sigma)$ is a subspace of the space of Kleinian groups $\mathcal K(\Sigma\times \mathbb R)$. The latter space is called the space of {\em quasi-Fuchsian} groups of the surface $\Sigma$ and is denoted as $\mathcal Q(\Sigma)$.  In other words quasi-Fuchsian groups are deformations of Fuchsian ones within the class of discrete subgroups of $PSL(2,\mathbb C)$. The space of quasi-Fuchsian groups is a subspace of the moduli space $\mathcal M(\Sigma,PSL(2,\mathbb C))$ of flat $PSL(2,\mathbb C)$ connections on the surface $\Sigma$.  If the surface $\sigma$ is not connected, the space $\mathcal Q(\Sigma)$ is by definition a product of the spaces of quasi-Fuchsian groups over the connected components of $\Sigma$. 

A Kleinian group acts not only on the hyperbolic space but also on its absolute boundary --- the complex projective space $\mathbb CP^1$. Its action is not properly discontinuous and the quotient is not a Hausdorff space. However if we remove the closure of the set of stable points of all elements of the group then the quotient becomes a Riemann surface homeomorphic to the ideal boundary $\partial X$ of $X$, thus establishing a map from the space of Kleinian groups $\mathcal K(X)$ to the Teichm\"uller space $\mathcal T(\partial X)$. The celebrated Bers uniformization theorem tells us that this map is an isomorphism unless the space $\mathcal K(X)$ is empty. In particular for quasi-Fuchsian groups the surface $\partial(\Sigma\times\mathbb R)$  is the double of the surface $\Sigma$ and thus there exist two projections from $\mathcal Q(\Sigma)$ to $\mathcal T(\Sigma)$. Once we fix the orientation of $\Sigma$ we can distinguish between these two projections. 

The space $\mathcal Q(\Sigma)$ is an open subset of the moduli space of $PSL(2,\mathbb R)$ flat connection on $\Sigma$ and hence inherits canonical complex symplectic structure.  It is also an open subset (in two ways) of the moduli space of complex projective structures on $\Sigma$ and hence the fibers of both projections to $\mathcal T(\Sigma)$ are Lagrangian \cite{Goldman}.

Consider  a manifold $X$ such that for every connected component of $\partial X$ the fundamental groups of it imbeds into the fundamental group of $X$ (if this condition is satisfied we say that the boundary is {\em incompressible}). A subgroup of a Kleinian group of $X$ corresponding to every boundary component is obviously quasi-Fuchsian. A collection of quasi-Fuchsian groups one for each boundary component of $\partial X$ defines a map $\mathcal L : \mathcal K(X) \rightarrow \mathcal Q(\partial X)$. Of course the map of $\mathcal K(X)$ to $\mathcal T(\partial X)$ defined above is just the composition of themap $\mathcal L$ with the natural projection $\mathcal Q(\partial X) \rightarrow \mathcal T(\partial X)$ defined by the orientation of $\partial X$ induced by the orientation of $X$. The Bers theorem states that the image of $\mathcal L$ is a section of this projection. Moreover one can show that this section is Lagrangian.

Unlike Teichm\"uller spaces which are quite well describable in coordinates and quantizable, our knowledge of the structure of the spaces of quasi-Fuchsian groups is rather scarce. In this paper we try to show that the question of description of the space of quasi-Fuchsian groups is related to affine Lie algebra $\widehat{PSL(2,\mathbb R)}$ like the description of ordinary Teichm\"uller space is related to the finite-dimensional one $PSL(2,\mathbb R)$. The main tool to study this question is the cosh-Gordon equation on a Riemann surface generalizing the Liouville equation.

\paragraph{Embedded surfaces.} Now recall a few facts about 2D surfaces embedded into a 3D Riemann manifold. The restriction of the 3D metric to the surface gives the metric on the surface, which is often called the {\em first fundamental form}. Consider a family of embeddings obtained from the original one by shifting each point along the geodesic orthogonal  to the image to the distance $\epsilon$ in the direction prescribed by the orientation of the surface. The first derivative of the induced metric w.r.t. $\epsilon$ at $\epsilon=0$ is called the {\em second fundamental form} of the embedded surface. The condition that the surface is minimal amounts to the requirement that the ratio of the second and the first quadratic forms is traceless. If our 3D manifold has constant curvature then the two fundamental forms satisfy further restrictions called the Gauss-Codazzi equation. The particular form of this equation of course depends on the sign of the curvature. 

\paragraph{Phase space.} Now let us come back to our original question about the phase space of 3D gravity. It is well known that the space of solutions of  Einstein equations in 3D with Euclidean signature and cosmological constant -1 is just the space of metrics with curvature -1. Solutions of the Einstein equateions are extrema of the Hilbert action
$$S(G)=\int (R(G)-1)d\mathrm{Vol},
$$ 
where $R(G)$ is the scalar curvature of the metric $G$ and $d\mathrm{Vol}=\sqrt{\det G}d^3x$ is the volume form induced by the metric $G$.

By definition the phase space of the theory corresponding to a 2D surface $\Sigma$ is the limit of the space of solutions of Einstein equations in a tubular neighborhood of the surface when we shrink the neighborhood. This space is considered modulo the kernel of the 2-form induced by the Hilbert action.  Since Einstein equations are of the second order its solution in a sufficiently small vicinity of $\Sigma$ is defined once we know it on the first infinitesimal neighborhood of $\Sigma$. The sufficiently small tubular neighborhood is fibered by geodesics orthogonal to the surface and, thus, is isomorphic to a direct product $\Sigma \times [-\varepsilon, \varepsilon]$ with the metric $G=g+tb+dt^2 + O(t^2)$. Here $g$ is a metric on $\Sigma$, $b$ is a quadratic form on $\Sigma$ and $t$ is the natural coordinate on the second factor. The forms $g$ and $b$ are obviously just the first and the second fundamental forms, respectively. A pair $(g,b)$ of forms defines a point in our phase space.

Recall that any Lagrangian induces an exact 2-form on the corresponding phase space as a differential of the boundary terms of its first variation. 
The canonical 2-form coming from the Hilbert action can be computed explicitly in terms of $g$ and $b$ as
$$
\omega = \delta \int_\Sigma (g^{ij}\delta b_{ij}+b^{ij}\delta g_{ij})d\mathrm{Vol}.
$$

The integral submanifolds of the kernel of the symplectic form $\omega$ are classes of an equivalence relation which can be described as follows. Two pairs $(g_1,b_1)$ and $(g_2,b_2)$ are considered as equivalent if there exists a 3D manifold $X$ with curvature -1 metric and two homotopic embeddings $\varphi_{1,2}:\Sigma \to X$, such that $(g_1,b_1)$ and $(g_2,b_2)$ are the pairs of fundamental forms induced by the embeddings $\varphi_1$ and $\varphi_2$, respectively. The corresponding equivalence relation can be considered as a groupoid action, and the kernel of the 2-form is generated by the corresponding algebroid.

In order to describe the quotient by this equivalence relation it is sufficient to pick one point in every orbit. We suggest to require the minimality of the surface $\Sigma$ in its first infinitesimal neighborhood. In the language of fundamental forms it amounts to the requirement that the ratio of the second and the first fundamental forms is traceless: $g^{ij}b_{ji}=0$.

This condition is justified by the fact that any complete hyperbolic structure on $\Sigma \times \mathbb R$ admits a minimal section (cf. \cite{Hass}). Unfortunately though such minimal sections are always isolated they are not always unique. We suggest a condition to single out one of them, but recognize that this condition is not so much justified by anything but that it looks nicely.

Obviously once we have a metric (the first fundamental form) on the surface it defines a complex structure on it. The minimality condition implies that the second fundamental form must be a sum of complex conjugated differentials $t$ and $\bar t$ of type (2,0) and (0, 2), respectively. We shall show that the Gauss-Codazzi condition amounts to holomorphicity of the quadratic differential $t$ as well as to a relation between $t$ and the conformal factor $e^\phi$ of the first fundamental form called the cosh-Gordon equation.

\paragraph{Hitchin hyper-K\"ahler structures.} The celebrated result of N.Hitchin \cite{Hitchin} tells us that given a compact simple Lie group $G$ for almost any flat $G^{\mathbb C}$ connection $\mathcal A$ on a Riemann surface with a {\em complex} simple Lie group $G^{\mathbb C}$ one can find a unique family of connections
\begin{equation}\label{Hitchin}
\mathcal A(\lambda) = \lambda^{-1}\Phi + A +\lambda \Phi^*
\end{equation}
being a Laurent polynomial of the parameter $\lambda \in \mathcal C$ and such that
\begin{enumerate}
 \item $\mathcal A(1)$ is $G$-gauge equivalent to the initial connection $\mathcal A$.
 \item The first and the last coefficients are of type (1,0) and (0,1), respectively.
 \item $\mathcal A(\lambda)=\mathcal A(-1/\bar{\lambda})^*$. (Here $*$ is the Hermitian involution.)
\end{enumerate}

The space of such families can be mapped isomorphically to the cotangent bundle to the space $Bun(\Sigma,G^{\mathbb C})$ of holomorphic $G^\mathbb{C}$-bundles on the surface $\Sigma$ (see e.g. \cite{Hitchin}). (The space $Bun(\Sigma,G^{\mathbb C})$ is isomorphic to the space $\mathcal M(\Sigma,G)$ of the moduli of $G$-connections on $\Sigma$ by the Narasimhan-Seshadri theorem).  Indeed the (0,1) part of the coefficient $A$ of the connection $A(\lambda)$ defines a holomorphic structure on the bundle, however the coefficient $\Phi$ defines a contangent vector to the space of holomorphic bundles.

This connection gives a family of maps from the moduli $\mathcal M(\Sigma,G^\mathbb C)$ of complex flat connections to itself depending on the parameter $\lambda \in \mathbb C^*$. Hitchin has shown that the family of complex and symplectic structures induced on $\mathcal M(\Sigma, G^\mathbb C)$ by this family is a hyper-K\"ahler one. In the limit $\lambda \rightarrow 0,\infty$ these structures tend to the ones on the cotangent bundle to the moduli $\mathcal M(\Sigma,G)$ of $G$-connections on $\Sigma$ with compact Lie group $G$. Such a family of connections can be considered as a single connection with values in an affine group $\widehat G^\mathbb C$.

\paragraph{Toda integrable system.} Recall that given an $r\times r$ Cartan matrix $C_{\alpha\beta}$ of a simple or affine Lie algebra $\mathfrak g$, the corresponding Toda equation is the following equation for an $r$-tuple of functions $\phi_1,\ldots,\phi_r$ on a complex plane:
$$
\frac{\partial^2 \phi_\alpha}{\partial z\partial \bar{z}}=\sum_\beta C_{\alpha\beta}e^{\phi_\beta}
$$
In particular for $\mathfrak g=\mathfrak{sl}(2)$ one gets the Liouville equation $\frac{1}{2}\frac{\partial^2 \phi}{\partial z\partial \bar{z}}=e^\phi$ and for the affine Lie algebra $\widehat{\mathfrak{sl}(2)}$ one gets the sinh-Gordon equation $\frac{1}{2}\frac{\partial^2 \phi}{\partial z\partial \bar{z}}=e^\phi-e^{-\phi}$.

There exists a straightforward way to reduce solution of these equations to the problem of Gauss factorization in the corresponding Lie group $G$ (cf. \cite{Leznov}). By Gauss factorization we mean the decomposition of an element $U$ of the group $G$ into a product $U_+U_-$, where $U_+$ and $U_-$ are elements of the upper and lower triangular Borel subgroups, respectively. Gauss decomposition exists for $U$ sufficiently close to the identity and is unique up to the transformation $(U_+,U_-)\mapsto(U_+K^{-1},KU_-)$, where $K$ is an element of the Cartan subgroup of the group $G$. 

One can easily check that the Toda system is equivalent to the flatness of the connection $A=(\sum_\alpha a^\alpha e_\alpha + \sum_\alpha b^\alpha h_\alpha)dz+(\sum_\alpha \bar{a}^\alpha f_\alpha + \sum_\alpha b^\alpha h_\alpha)d\bar{z}$, where $e^{\phi_\alpha}=a^\alpha \bar a^\alpha$ and $\{e_\alpha,h_\alpha,f_\alpha\}$ are the standard Cartan-Weyl generators of the Lie algebra $\mathfrak g$. Let $A^0=(\sum_\alpha a^\alpha e_\alpha + \sum_\alpha b^\alpha h_\alpha)dz$ be any holomorphic connection taking value in a subspace of $\mathfrak g$ generated by the Cartan subalgebra and positive simple roots. To any such connection we are going to associate a solution of the Toda system in the following way. Let $g(z)=\mathrm{Pexp}\int_0^z A^0$ be the holomorphic group-valued function satisfying $g(0)=1$ and $dg = A^0g$. Let $\bar{A}^0=(\sum_\alpha \bar{a}^\alpha f_\alpha + \sum_\alpha \bar{b}^\alpha h_\alpha)d\bar{z}$ and $g^*(\bar{z})$ be an antiholomorphic solution of $g^*(0)=1$ and $dg^* = \bar{A}^0g^*$. Let finally $\mathbf g(z,\bar{z}) = (g^{-1}g^*)_-g$. Then the connection $A=d\mathbf g \mathbf g^{-1}$ is the connection providing a solution to the Toda system.

\section{Cosh-Gordon metrics.}

Let $\Sigma$ be a closed two dimensional surface of genus $g$ and $\mathcal T(\Sigma)$ be the corresponding Teichm\"uller space. Recall that the cotangent space to $\mathcal T(\Sigma)$ at a point $S$ can be canonically identified with the space of holomorphic quadratic differentials on $S$. 

Let $(S,t)\in T^*\mathcal T(\Sigma)$ be a pair of a Riemann surface $S$ and a holomorphic quadratic differential on $S$. We shall say, that an Hermitian metric $\rho$ is a {\em cosh-Gordon} one if it satisfies the relation:

\begin{equation}\label{CG}
\frac{1}{2}\partial\bar{\partial}\ln \rho = \rho + t\bar{t}\rho^{-1},
\end{equation}
and the condition
\begin{equation}\label{ineq}
 \rho > t\bar{t}\rho^{-1}
\end{equation}

These relations have the following features:

1. Though the relations are written in local coordinates, the dimension of all terms are equal to (1,1), and thus this condition is coordinate independent since the l.h.s. is just minus the curvature of the metric $\rho$ times the volume form.

2. If the quadratic differential $t$ is equal to zero, the relation becomes the Liouville equation, which means just that the metric $\rho$ has curvature $-1$.

3. In any neighborhood without zeroes of $t$ on can choose a local coordinate function $z=\int\sqrt{t}$ making $t=1$. The substitution $\rho = e^\phi$ turns (\ref{CG}) into the usual cosh-Gordon equation. 

4. Integrating both sides of the cosh-Gordon equation (\ref{CG}) over the whole surface, using the Gauss-Bonnet theorem for the integral of the l.h.s. and Cauchy-Schwartz inequality for the r.h.s. one obtains the inequality:
\begin{equation}
-\pi\chi(\Sigma) \geq \int_\Sigma\sqrt{t\bar{t}}, 
\end{equation}
where  $\chi(\Sigma)$ is the Euler characteristic of $\Sigma$. Thus a
cosh-Gordon metric can exist only if the quadratic differential is not
too large and only if the Euler characteristic $\chi(\Sigma)$ is negative.

5. Let $\rho$ be a solution of (\ref{CG}). Taking $\tilde{\rho} = \rho e^\phi$ for a function $\phi$, substituting this expression into  (\ref{CG}) and linearizing one gets
\begin{equation}\label{linear}
\partial\bar{\partial}\phi = (\rho-t\bar{t}\rho^{-1})\phi.
\end{equation}
This equation has no nonvanishing solutions provided the inequality (\ref{ineq}) is satisfied. Thus any cosh-Gordon metric is isolated.  (Caution:  a solution of \ref{CG} may be isolated even if the inequality (\ref{ineq}) is not satisfied.)

6. Let $G$ be a metric on $\Sigma \times \mathbb R = \{(z,l)\}$ given by:
\begin{equation}\label{metric}
G = dl^2 + \rho\cosh^2 l -t\bar{t}\rho^{-1}\sinh^2 l +  (t + \bar{t})\cosh l\sinh l
\end{equation}
Then 

$\rho$ is cosh-Gordon if and only if the metric $G$ is of curvature -1. In this case we also have that the section $l=0$ is a minimal area surface with respect to the metric $G$ and the lines corresponding to fixed $z$ are geodesics orthogonal to the section $l=0$.

In order to prove these properties let us first prove the following
\begin{lemma}
Let $A$ be a flat connection  1-form on a surface $\Sigma$ taking values in the Lie algebra $\mathfrak{sl}(2,\mathbb C)$ in the standard two-dimensional representation and  $\Lambda(l) = \left( \begin{array}{cc}e^l& 0\\ 0& e^{-l}\end{array}\right)$. Let $G$ be the metric on $\Sigma \times \mathbb R=\{(z,l)\}$ given by
\begin{equation}\label{indmetric}
G(A)=\frac{1}{4}\det(\Lambda^{-1}A\Lambda + \Lambda^{-1}d\Lambda + A^+),
\end{equation}
Then this metric $G(A)$ has constant curvature -1 provided it is nondegenerate. Moreover the lines $(z,\cdot)$ are in this case geodesics for $G$. 

The metric $G(A)$ does not change if we change the connection $A$ by a gauge transformation 
\begin{equation}\label{gauge}
A \mapsto h^{-1}dh + h^{-1}A h,
\end{equation}
where $h$ is a map from $\Sigma$ to diagonal unitary matrices.
\end{lemma}

{\em Proof of the lemma.~~}  Identify the hyperbolic space   $\mathbb{H}^3$ with the space of Hermitean $2\times2$ matrices with $\det H =1$ and the metric given by  $\frac{1}{4}\det dH$ (here $dH$ is considered as  a function on the tangent bundle to $\mathbb{H}^3$). The group $PSL(2,\mathbb C)$ acts on  $\mathbb{H}^3$ by $(g,H)\mapsto gHg^+$. Here $g\in PSL(2,\mathbb C)$ and  $g^+$ is the Hermitean conjugate of $g$.  Consider a curve in $\mathbb H^3$ given by $\Lambda(l) = \left( \begin{array}{cc}e^l& 0\\ 0& e^{-l}\end{array}\right)$, where $l \in \mathbb R$. One can easily check that this curve is a geodesic and $l$ is a natural parameter on it. Any other geodesic is of course an image of this one under the action of $PSL(2,\mathbb{C})$.

Let $U$ be a domain in the surface $\Sigma$ and let now $g$ be a map  $g: U \rightarrow PSL(2,\mathbb{C})$ such that $A=g^{-1}dg$. Such a map $g$ induces a map $U \times \mathbb{R} \rightarrow \mathbb{H}^3$ by $(z,l) \mapsto g(z) \Lambda(l) g^+(z)$.

The induced metric on $U \times \mathbb{R}$ reads as 
$$ G =\det(d( g\Lambda g^+)) =
\det(\Lambda^{-1}g^{-1}dg \Lambda + \Lambda^{-1}d\Lambda + dg^+
(g^+)^{-1}). $$
which coincides with (\ref{indmetric}). It is by construction of curvature -1 provided the map $U \times \mathbb{R} \rightarrow \mathbb{H}^3$ is an embedding, what follows from the nondegeneracy of $G$.

The image of a real line $(z,\cdot)$, where  $z \in U$ is a geodesic in  $\mathbb{H}^3$ since $\Lambda(\cdot)$ is a geodesic.

This map and therefore the metric is obviously invariant w.r.t. the transformation $g \mapsto gh$, where $h$ is a diagonal unitary
matrix: $ h =  \left( \begin{array}{cc}e^{i\psi}&0\\
0&e^{-i\psi}\end{array}\right)$, and $\psi \in \mathbb{R}$. 

The Lemma 1 is proven.

In order to prove the property 6 of the cosh-Gordon mertic relating cosh-Gordon equation on $\Sigma$ with the condition on the curvature of the 3D meric $G$ it is now sufficient to find a flat connection $A$ determined by a solution $\rho$ of the cosh-Gordon equation and show that this connection gives  the metric .  Here is the formula for this connection. 

\begin{equation}\label{CGconnection}
A =
 \left( \begin{array}{cc}-\frac{1}{4}\partial \phi& te^{-\phi/2} \\
e^{\phi/2}&\frac{1}{4}\partial\phi\end{array}\right)dz +
 \left( \begin{array}{cc}\frac{1}{4}\bar{\partial} \phi&-te^{-\phi/2} \\
e^{\phi/2}&-\frac{1}{4}\bar{\partial}\phi\end{array}\right)d\bar{z},
\end{equation}
where $\phi = \log \rho$.

The verification that this connection is indeed flat and that it gives the metric (\ref{metric}) is  a straightforward substitution of ( \ref{CGconnection}) into (\ref{indmetric}).

To show where the connection (\ref{CGconnection}) came from let us prove the following
\begin{lemma}
The connection $A$ given by the expression (\ref{CGconnection}) with is singled out by the following three conditions on an $A$ and the metric $G$ induced on $\Sigma\times \mathbb R$ via (\ref{indmetric}):
\begin{enumerate}
\item $A$ is flat.
\item The section $l=0$ is minimal.
\item The lines $z=z_0$ are orthogonal to the section $l=0$.
\end{enumerate}
\end{lemma}

In order to prove the Lemma let us look for  the desired connection in the form
$$\mathcal{A} = \left( \begin{array}{cc}\xi&\zeta\\
\bar{\eta}&-\xi\end{array}\right),$$ where $\xi,\zeta$ and $\eta$ are one-forms.

Substituting this expression for the connection $A$ into (\ref{indmetric}) we get the expression for the 3D metric:
$$G = (dl + \xi +\bar{\xi})^2 + (e^{-l}\zeta +
e^l\eta)(e^{-l}\bar{\zeta} + e^l\bar{\eta}).
$$
The $U(1)$ gauge transformation \ref{gauge} reads as:
$$ (\mathcal{A},\psi) \mapsto h^{-1}dh + h^{-1}\mathcal{A}h
=\left( \begin{array}{cc}\xi+id\psi&\zeta e^{2i\psi}\\ \bar{\eta} e^{-2i\psi}
&-\xi-id\psi\end{array}\right) $$

Introduce new variables $\alpha = (\zeta + \eta)/2$ and  $\beta = (\zeta
- \eta)/2$.  In terms of these variables the metric takes the form
$$ G =
\alpha\bar{\alpha}\cosh^2 l + \beta\bar{\beta}\sinh^2 l -
(\alpha\bar{\beta}+\beta\bar{\alpha})\cosh l \sinh l + (dl + \xi
+\bar{\xi})^2.
$$

The orthogonality condition \textit{3} is obviously just
\begin{equation} \label{ortho}
\xi + \bar{\xi} = 0
\end{equation}

The minimality condition \textit{2} is
\begin{equation}\label{mini}
\alpha\wedge\bar{\beta}+\beta\wedge \bar{\alpha}=0.
\end{equation}

Substituting these conditions into the zero curvature equation $dA + A \wedge A=0$ on gets that $\alpha\wedge\bar{\beta}=0$ and hence the forms $\alpha$ and $\bar{\beta}$ are proportional.

Under these conditions the zero curvature  equation reduces to the three equations:
$$
  \begin{array}{lll}
          d\alpha &=& 2\alpha \wedge \xi,\\
          d\beta &=& 2\beta \wedge \xi,\\
          d\xi &=& (\beta\wedge \bar{\beta} - \alpha\wedge \bar{\alpha}).
  \end{array}
$$

Introduce a complex coordinate $z$ such that the form $\alpha$ is proportional to $dz$. Using the gauge freedom one can always make the proportionality coefficient to be real, and it must be nonzero since the metric $G$ is nondegenerate.  Therefore we can write the form $\alpha$ as $e^{\phi/2}dz$, where $\phi$ is a real-valued function, and $\beta$ as $\bar{t}e^{-\phi/2}d\bar{z}$, where $t$ is a complex-valued function. Substituting this Ansatz into the flatness condition one easily gets that 
$$\xi = \frac{1}{4}(\partial\phi dz- \bar{\partial}\phi d\bar{z})$$
$$\bar{\partial}t = 0$$
\begin{equation}
\frac{1}{2}\partial\bar{\partial}\phi = e^\phi + t\bar{t}e^{-\phi}.
\end{equation}
So substituting the expressions for $\xi,\eta, \zeta$ one obtains the desired form (\ref{CGconnection}) of the connection $A$.

\section{Flat connection with affine gauge group.}
It turns out that the connection (\ref{CGconnection}) can be included into a family of flat connections depending on a parameter $\lambda \in \mathbb C^*$ analogous to the Hitchin family (\ref{Hitchin}). However in our construction the complex structure on the surface is a function of the original connection, and is not fixed once for all connections as in Hitchin's construction. Another difference is that our families are related to pseudohermitean involution on $PSL(2,\mathbb C)$, while Hitchin's families use the Hermitean one. The family also can be interpreted as a connection in a bundle with affine group $\widehat{ SL(2,\mathbb C)}$ as a structure group. It turns out that the connections we constructed coincide with the ones used to describe local solutions of affine Toda equations discussed in the section 1.

\begin{lemma}
Let $\mathcal A(\lambda)$ be a family of $PSL(2,\mathbb C)$-connections depending on a parameter $\lambda$ and satisfying the following conditions:
\begin{enumerate}
\item Flatness: $dA(\lambda)+A(\lambda)\wedge A(\lambda)=0$ for any value of $\lambda \in \mathbb C^*$.
\item Polynomiality: $A(\lambda)=\Phi\lambda^{-1}+A+\bar{\Phi}\lambda$ is a Laurent polynomial in $\lambda$ of degree $[-1,1]$.
\item Degeneracy: At least one component of the 1-form $\Phi$ is a degenerate matrix.
\item Reality: $\mathcal{A}(-\bar{\lambda}^{-1})^* =
-\mathcal{A}(\lambda)$, where $^*$ is the pseudohermitean conjugation $g^* =
Cg^+C$ and $C =   \left( \begin{array}{cc}1&0\\0&-1\end{array}\right)$.
\end{enumerate}
Then for an open subset of such connections one can choose a gauge and a complex structure on the surface, and a local complex coordinate $z$ such that the connection takes the form 

\begin{equation}\label{matrix-connection}
\mathcal{A}(\lambda) =
 \left( \begin{array}{cc}-\frac{1}{4}\partial \phi& \lambda te^{-\phi/2} \\
e^{\phi/2}&\frac{1}{4}\partial\phi\end{array}\right)dz +
 \left( \begin{array}{cc}\frac{1}{4}\bar{\partial} \phi& e^{\phi/2} \\
-\lambda^{-1}\bar{t}e^{-\phi/2}&-\frac{1}{4}\bar{\partial}\phi\end{array}\right)d\bar{z}.
\end{equation}
where $t$ is a holomorphic quadratic differential and $\rho=e^\phi$ is a solution for the cosh-Gordon equation (\ref{CG}).
\end{lemma} 

{\em Proof of the lemma.~~}  The $\lambda^{-2}$-term of the flatness condition amounts to $\Phi\wedge\Phi=0$ what implies that the components $\Phi_1$ and $\Phi_2$ of the 1-form $\Phi$ commute, and since they belong to a rank 1 algebra, they are proportional. (The reality condition implies that $\bar{\Phi}=\Phi^*$ and thus the components of $\bar{\Phi}$ are proportional and degenerate as well.)  Thus one can choose a gauge making $\Phi$ proportional to $E=\left( \begin{array}{cc}0&1\\0&0\end{array}\right)$. The proportionality coefficient is a 1-form. For an open subset of connections this 1-form is nowhere real on $\Sigma$, and this will be the subset of all connections which we are going to consider. Introduce a complex coordinate $z$ on $\Sigma$ and write this form as $adz$. Thus in this gauge the connection $\mathcal A(\lambda)$ can be written as
\begin{equation}\label{long-connection}
\begin{array}{rl}
\mathcal A(\lambda)&=\left((e+\lambda^{-1} a)E+hH+fF\right)dz+\\ &+\left((\bar{e}-\lambda\bar{a})F-\bar{h}H+\bar{f}E\right)d\bar{z},
\end{array}
\end{equation}
where $a,b,e,f,h$ are complex functions on the surface and $E,F,H$ are the standard Cartan-Weyl generators of the Lie algebra $\mathfrak{sl}(2)$, satisfying $[E,F]=H$, $ [H,E]=2E$, $[H,F]=-2F$, $E^*=-F$, $F^*=-E$, $H^*=H$. 

The zero curvature curvature condition for of this connection reads as
$$0= d\mathcal A(\lambda)+\mathcal A(\lambda)\wedge \mathcal A(\lambda)=
$$ 
$$=((e+\lambda^{-1}a)(\bar{e}-\lambda\bar{a})-f\bar{f}-\bar{\partial} h-\partial\bar{h})Hdz\wedge d\bar{z}+
$$
$$+(2h\bar{f}+2(e+\lambda^{-1} a)\bar{h}+\partial\bar f+\bar\partial e-\lambda^{-1}\bar{\partial}a)Edz\wedge d\bar{z}+
$$
$$+(-2f\bar{h}-2h(\bar{e}-\lambda\bar{a})+\partial\bar e-\lambda\partial\bar a-\bar{\partial}f)Fdz\wedge d\bar{z}.
$$

Equating to zero the terms at $E,F,H$ and with different powers of the parameter $\lambda$ one gets a system of equations. The terms at $\lambda H$ and $\lambda^{-1}H$ give $a\bar{e}=\bar{a}e=0$. In an open subset of connections $a\neq 0\neq \bar{a}$, and thus $e=\bar{e}=0$. Taking into account these relations one gets
\begin{eqnarray*}
\bar\partial a&=&2\bar ha\\
\partial\bar a&=&2h\bar a\\
\bar\partial f&=&-2\bar hf\\
\partial\bar f&=&-2h\bar f\\
\partial\bar h +\bar\partial h &=& -(a\bar a+f\bar f)
\end{eqnarray*}

Restrict further the set of connection by the requirement $f\neq 0$. One can use a gauge transformation by a diagonal matrix in order to make it real and positive: $f=\bar f\ge 0$. Using this choice introduce new complex function $t=af$ and a real one $\phi = 2\log f$ and transform the equations into  
\begin{eqnarray*}
\bar\partial t&=&0\\
e^\phi+t\bar{t}e^{-\phi}-\frac12\bar{\partial}\partial\phi&=&0,
\end{eqnarray*}
Substituting the expressions of $a,f,h$ in terms of $t$ and $\phi$ into (\ref{long-connection}) one gets the desired connection (\ref{matrix-connection}).

\section{Sinh-Gordon metric and the phase space for signature 2+1.}

\paragraph{Sinh-Gordon metric.} Instead of cosh-Gordon metric one can consider a {\em sinh-Gordon} one, satisfying the equation
\begin{equation}\label{SG}
\frac{1}{2}\partial\bar{\partial}\ln \rho = \rho - t\bar{t}\rho^{-1},
\end{equation}
without any extra conditions.

Analogously to the cosh-Gordon equation (\ref{CG}) this equation defines a Hermitean metric with curvature $-1$ if $t=0$. It transforms into the ordinary sinh-Gordon equation for a coordinate system where $t=1$. The linearization of this equation

$$
\partial\bar{\partial}\phi = (\rho+t\bar{t}\rho^{-1})\phi
$$
has no hontrivial solution and therefore any sinh-Gordon metric is isolated.

The geometric meaning of the sinh-Gordon equation is characterized by the following
\begin{lemma}.

1. The sinh-Gordon equation is equivalent to the requirement for the metric $\rho+t+\bar{t}$ to have curvature $-1$.

2. The map sending $t$ to the Riemann surface with complex structure defined by $\rho+t+\bar{t}$ is an isomorphism between $T^*_S\mathcal T(\Sigma)$ and $\mathcal T(\Sigma)\times \mathcal T(\Sigma)$.
\end{lemma}

This proposition is a part of common knowledge and we do not give a detailed proof of it here. The statement \textit{1} can be proven by direct computation of the curvature. The statement \textit{2} splits into two statements: One claims that for any $(S,t)$ in $T^*\mathcal T(\Sigma)$ there exists a sinh-Gordon metric. Another one claims that for any $S_1,S_2\in \mathcal T(\Sigma)$ there exists a quadratic differential $t$ on $S_1$ and sinh-Gordon metric $\rho$, such that the metric $\rho+t+\bar{t}$ defines the complex structure coinciding with the given one on $S_2$. 

The first of these sub-statements follows from standard variational arguments, since the action functional, corresponding to the sinh-Gordon equation is bounded below. The second one follows from the existence and uniqueness theorem for harmonic maps between spaces of negative curvature. Indeed, let $\rho_1$ be a Hermitian metric on $S_1$ of curvature -1. Then a map $h:S_2\rightarrow S_1$ is harmonic iff the (2,0) part $t$ of the metric $h^* \rho_1$ called {\em Hopf differential} is holomorphic. Thus the (1,1) part of $h^* \rho_1$ is sinh-Gordon w.r.t. the Hopf differential $t$. 

The other properties of sinh-Gordon metrics are similar to those of cosh-Gordon ones:

1. The signature 2+1 metric on $\Sigma \times \mathbb R$ given by
\begin{equation}\label{sg-metric}
G = -dl^2 + \rho\cos^2 l + t\bar{t}\rho^{-1}\sin^2 l +  (t + \bar{t})\cos l\sin l
\end{equation}
has curvature -1 if and only if the equation (\ref{SG}) is satisfied. In this case the surface $l=0$ is minimal and the lines $(z,\cdot)$ are geodesics orthogonal to it.

2. The equation (\ref{SG}) is equivalent to the flatness of the connection $A$ given by
\begin{equation}\label{sg1-matrix-connection}
\mathcal{A}_1(\lambda) =
 \left( \begin{array}{cc}-\frac{1}{4}\partial \phi& -i\lambda te^{-\phi/2} \\
e^{\phi/2}&\frac{1}{4}\partial\phi\end{array}\right)dz +
 \left( \begin{array}{cc}\frac{1}{4}\bar{\partial} \phi& e^{\phi/2} \\
i\lambda^{-1}\bar{t}e^{-\phi/2}&-\frac{1}{4}\bar{\partial}\phi\end{array}\right)d\bar{z}.
\end{equation}

or, equivalently, of the connection
\begin{equation}\label{sg2-matrix-connection}
\mathcal{A}_2(\lambda) =
 \left( \begin{array}{cc}\frac{1}{4}\partial \phi& -i\lambda te^{-\phi/2} \\
-e^{\phi/2}&-\frac{1}{4}\partial\phi\end{array}\right)dz +
 \left( \begin{array}{cc}-\frac{1}{4}\bar{\partial} \phi& -e^{\phi/2} \\
i\lambda^{-1}\bar{t}e^{-\phi/2}&\frac{1}{4}\bar{\partial}\phi\end{array}\right)d\bar{z}.
\end{equation}

Observe that these two connections are $SU(1,1)$-connections for $\lambda=1$ with independent monodromies.

\end{document}